\documentclass[12pt,reqno]{amsart}
\usepackage{hyperref,cleveref,amssymb,eqnarray,mathrsfs,xcolor}

\theoremstyle{plain}
\newtheorem{theorem}{Theorem}[section]
\newtheorem{proposition}[theorem]{Proposition}
\newtheorem{corollary}[theorem]{Corollary}
\newtheorem{lemma}[theorem]{Lemma}

\theoremstyle{definition}

\newtheorem{question}[theorem]{Question}
\newtheorem{example}[theorem]{Example}

\newcommand{\Pelc}[1]{\left(\Sigma X_n\right)_p}
\DeclareSymbolFont{bbold}{U}{bbold}{m}{n}
\DeclareSymbolFontAlphabet{\mathbbold}{bbold}

\DeclareMathOperator{\ran}{ran}

\DeclareMathOperator{\diag}{Diag}

\renewcommand{\le}{\leqslant}

\newcommand{\cN}{\mathcal{N}}
\newcommand{\cR}{\mathcal{R}}

\begin{document}

\title[Positive self-commutators]{Positive self-commutators of positive operators}%\\

\author{Roman Drnov\v sek}
\address{Faculty of Mathematics and Physics, University of Ljubljana,
  Jadranska 19, 1000 Ljubljana, Slovenia \ \ \ and \ \ \
  Institute of Mathematics, Physics, and Mechanics, Jadranska 19, 1000 Ljubljana, Slovenia}
\email{roman.drnovsek@fmf.uni-lj.si}

\author{Marko Kandi\'c}
\address{Faculty of Mathematics and Physics, University of Ljubljana,
  Jadranska 19, 1000 Ljubljana, Slovenia \ \ \ and \ \ \
  Institute of Mathematics, Physics, and Mechanics, Jadranska 19, 1000 Ljubljana, Slovenia}
\email{marko.kandic@fmf.uni-lj.si}


\date{\today}

\begin{abstract}
We consider a positive operator $A$ on a Hilbert lattice such that its self-commutator $C = A^* A - A A^*$ is positive. If $A$ is also idempotent, then it is an orthogonal projection, and so $C = 0$.
Similarly, if $A$ is power compact, then $C = 0$ as well. We prove that every positive compact central operator on a separable infinite-dimensional Hilbert lattice $\mathcal H$ is a self-commutator of a positive operator.
We also show that every positive central operator on $\mathcal H$ is a sum of two positive self-commutators of positive operators.
\end{abstract}

\maketitle

\section{Introduction and Preliminaries}\label{Introduction}

Positive commutators of positive operators on Banach lattices have been the subject of extensive research across various contexts. A systematic investigation into their
properties began with \cite{Bracic:10}, where the authors examined the spectral properties of
the positive commutator $[A, B] := AB - BA$ formed by positive compact operators $A$ and $B$. They established that such a commutator is always quasinilpotent, and it is contained in the radical of the Banach algebra generated by 
$A$ and $B$.

Further developments were seen in \cite{DK11}, where positive commutators of either positive nilpotent or positive power compact operators were explored. Independently, Drnov\v sek \cite{Drnovsek:11} and Gao \cite{Gao:14} proved that a positive commutator of positive operators, provided one of them is compact, is necessarily quasinilpotent. Kandi\'c and \v Sivic proved in \cite{KS17} that such commutator necessarily belongs to the radical of the Banach algebra generated by $A$ and $B$. In \cite{KS17b} they were investigating the dimension of the algebra generated by two positive $n\times n$ matrices. The most interesting result yields that, under some technical conditions, the unital algebra generated by two positive idempotent operators with a positive commutator is at most nine dimensional. Drnovšek obtained the same result in \cite{Drnovsek:18} 
under different assumptions. 
% He proved that two positive idempotent operators with a positive commutator on a vector lattice whose order dual separates the points is again at most nine dimensional. 
In \cite{DK19} the authors proved that every positive compact operator on a separable Banach lattice $L^p(\mu)$ ($1\leq p<\infty)$ is a commutator of positive operators. Lastly, in \cite{DK25}  authors studied commutators greater than a perturbation of the identity in ordered normed algebras. 

In this paper, we explore positive self-commutators of positive operators on Hilbert lattices. Our results are motivated by Radjavi's work \cite{Radjavi:66}, which essentially proves that a self-adjoint operator $A$ on a separable Hilbert space is a self-commutator if and only if zero lies within the convex hull of the essential spectrum of $A$. Notably, Sourour \cite{Sourour:79} provided a concise proof of this result. 
According to \cite[Theorem 11.5]{Con00}, for every normal operator $N$ on a Hilbert space, there exists a measure space $(\Omega, \Sigma, \mu)$ and a function $\phi \in L^\infty(\mu)$ such that $N$ is unitarily equivalent to the multiplication operator $M_\phi$ acting on $L^2(\mu)$. Furthermore, if $\mu$ is $\sigma$-finite, then by \cite[Example 2.67]{AB06}, multiplication operators of the form $M_\phi$, where $\phi \in L^\infty(\mu)$, are precisely orthomorphisms, or equivalently, central operators on the Banach lattice $L^p(\mu)$ for $1 \leq p < \infty$. 
Therefore, the main goal of the paper is to study which positive central operators are self-commutators of positive operators.

The paper is organized as follows. In \Cref{First Observations}, we establish that a positive idempotent operator on a Hilbert lattice with a positive self-commutator is an orthogonal projection. Additionally, we show that the zero operator is the only positive operator on a Hilbert lattice that is a self-commutator of a positive power compact operator.
In \Cref{Positive isometries}, we prove that every positive bounded diagonal operator on a Hilbert lattice $\mathcal{H} \in \{\ell^2, \ell^2_n\}$ factors through $L^2[0,1]$ as a self-product of the form $X^*X$ for some positive linear operator $X \colon \mathcal{H} \to L^2[0,1]$. Furthermore, there exists a positive self-adjoint operator $Y \colon L^2[0,1] \to L^2[0,1]$ such that $Y^2 = XX^*$. This result proves to be particularly useful in \Cref{Positive central self-commutators}, where we prove the order analog of \cite[Theorem 5]{Radjavi:66} stating that every positive compact central operator on a separable infinite-dimensional Hilbert space is a self-commutator of a positive operator.
Finally, in \Cref{Section 5}, we establish that every positive central operator on a separable Hilbert lattice can be expressed as a sum of two positive self-commutators of positive operators.
 
Before we proceed to the results, we briefly recall some facts about Banach and Hilbert lattices and operators acting on them.    
Let $L$ be a vector lattice with the positive cone $L^+$. The band
$$ S^d :=\{x\in L:\; |x|\wedge |y|=0  \textrm{ for all }y \in S\}$$
is called the {\it disjoint complement} of a set $S$ of $L$.
A band $B$ of $L$ is said to be a {\it projection band} if $L=B\oplus B^d$.
An operator $A$ on $L$ is called {\it positive} if it maps the cone $L^+$ into itself.
An operator $A$ on $L$ is called {\it negative} if the operator $-A$ is positive.

Let $A$ be a positive operator on a vector lattice $L$. The {\it null ideal} $\cN(A)$ is the ideal in $L$ defined by
$$ \cN(A) = \{ x \in L : A |x| = 0 \} .  $$
When $\cN(A) = \{0\}$, we say that the operator $A$ is {\it strictly positive}.
The {\it range ideal} $\cR(A)$ of $A$ is the ideal generated by the range of $A$, that is,
$$ \cR(A) = \{ y \in L : \,  \exists x \in L^+ \textrm{\ \ such that \ } |y| \le A x \} . $$
A positive operator $A$ on a vector lattice is an \emph{orthomorphism} if every band of $L$ is invariant under $A$. If there exists a positive real number $\lambda$ such that 
$|Ax|\leq \lambda |x|$ for every $x\in L$, then $A$ is called \emph{central}. By \cite[Theorem 3.29]{AA02}, an operator on a Banach lattice is central if and only if it is an orthomorphism. 
% An operator $A$ on $L$ is called {\it order continuous} if every net $\{x_\alpha\}$ order converging to zero is mapped to the net $\{A x_\alpha\}$ order converging to zero as well.

Let $L^p(\mu)$ $(1\leq p<\infty)$ be a separable Banach lattice. Then by \cite[Theorem 7.1]{Bohnenblust} the Banach lattice $L^p(\mu)$ is isometric and order isomorphic to one of the following Banach lattices $\ell_{n}^p, \ell^p$, $L^p[0,1]$, $\ell^p\oplus L^p[0,1]$ or $\ell_n^p\oplus L^p[0,1]$. If $\mathcal H$ is a Hilbert lattice, then by \cite[Theorem IV.6.7]{sch}
it is isometric and order  isomorphic to a Hilbert lattice of the form $L^2(\Omega,\Sigma,\mu)$, where $\Omega$ is a locally compact Hausdorff space and $\mu$ is a strictly positive Radon measure. Therefore, every band $B$ in a Hilbert lattice satisfies $B^d=B^\perp$.
If $A$ is a bounded operator on a Hilbert lattice, then by \cite[Proposition 1.1]{Radl}, the operator $A$ is positive if and only if the Hilbert space adjoint operator $A^*$ is positive.

For the terminology and details not explained here we refer the reader to  \cite{AA02} or \cite{AB06}.

\section{First observations on positive self-commutators}\label{First Observations}

Let $A$ be a positive operator on a Hilbert lattice $\mathcal H$ such that its self-commutator $C = A^* A - A A^*$ is also positive.
Since $A^* A \geq A A^*$, a simple induction gives that $(A^* A)^n \leq (A^*)^n A^n$ for each positive integer $n$. It follows that
$\|(A^* A)^n \| \leq \|(A^*)^n\| \| A^n \|$, and so $r(A^* A) \leq r(A^*) r(A) = r(A)^2$, where $r$ denotes the spectral radius.
In particular, if $r(A) = 0$ then $r(A^* A) = 0$. As $A^* A$ is self-adjoint, we have $A^* A = 0$, and so $A=0$ and $C=0$.
Therefore, we can assume that $r(A)>0$.

Let us first consider the following question: which positive operators on a Hilbert lattice $\mathcal H$
are self-commutators of positive idempotent operators?

\begin{proposition}\label{idempotent}
Let $A$ be a positive idempotent operator on a Hilbert lattice $\mathcal H$ such that its self-commutator $C = A^* A - A A^*$ is also positive.
Then $A$ is an orthogonal projection, and so $C = 0$.
\end{proposition}

\begin{proof}
Assume first that $\cN(A) = \{0\}$ and $\cR(A)^d = \{0\}$.
Since closed ideals are bands in $\mathcal H$, the range ideal  $\cR(A)$ is dense in $\mathcal H$. Therefore, it follows from
$A C A = A A^* A - A A^* A = 0$ that $C=0$ and the operator $A$ is normal.
Furthermore,
$$ (A - A^* A)^* (A - A^* A) = A^* A - A^* A -A^* A +A^* A A^* A = 0 , $$
so that $A - A^* A =0$. It follows that $A^* = A^* A = A$
completing the proof in this case.

Consider now the general case. Let us define the bands $H_1$, $H_2$, $H_3$ and $H_4$ by $H_1 = \cN(A) \cap \cR(A)^d$,
$H_2 = \cN(A) \cap \cR(A)^{d d}$, $H_3 = \cN(A)^d \cap \cR(A)^{d d}$ and $H_4 = \cN(A)^d \cap \cR(A)^d$.
% \textcolor{red}{Ali treba utemeljiti, da je tudi $\cR(A)^{dd}$ invarianten za $A$?}
With respect to the band decomposition $\mathcal H = H_1 \oplus H_2 \oplus H_3 \oplus H_4$, the idempotent $A$ has the form
$$ A=\left( \begin{matrix}
   0   &  0  & 0  & 0  \\
   0   &  0  & X  & Z  \\
   0   &  0  & E  & Y  \\
   0   &  0  & 0  &  0
\end{matrix} \right)  , $$
where $E$, $X$, $Y$ and $Z$ are positive operators on the appropriate bands.
It follows from $A^2 = A$ that $E^2 = E$, $X E = X$, $E Y = Y$ and $X Y = Z$.
We also have
$\cN(E) = \{0\}$, $\cN(Y) = \{0\}$, $\cR(E)^d = \{0\}$ and $\cR(X)^d = \{0\}$.
We now compute
$$ A^* A = \left( \begin{matrix}
   0   &  0  & 0  & 0  \\
   0   &  0  & 0  & 0  \\
   0   &  0  & X^* X +E^* E  & (X^* X +E^* E) Y  \\
   0   &  0  & Y^*(X^* X +E^* E)  &  Y^*(X^* X +E^* E) Y
\end{matrix} \right) $$
and
 $$ A A^* = \left( \begin{matrix}
   0   &  0  & 0  & 0  \\
   0   &  X X^* + X Y Y^* X^*  & X(I + Y Y^*) E^* & 0  \\
   0   &  E(I + Y Y^*) X^*  & E E^* + Y Y^*  &  0  \\
   0   &  0  & 0  &  0
\end{matrix} \right)  . $$

Comparing the $(2,2)$-block in $C = A^* A - A A^* \geq 0$, we obtain that
$X X^* = 0$, and so $X = 0$.
Comparing the $(3,3)$-block, we conclude that $E^* E \geq E E^*+ Y Y^*\geq E E^*$.
By the special case, $E$ is an orthogonal projection, and therefore $Y Y^* = 0$, so that $Y = 0$.
It follows that $A$ is an orthogonal projection as well. This completes the proof.
\end{proof}

Suppose that a positive matrix $C$ is a self-commutator of a positive matrix. Then, \cite[Theorem 2.1]{Bracic:10} yields that $C$ is nilpotent. Since $C$ is also hermitian, it needs to be zero. This result can be extended to positive power compact operators on Hilbert lattices.
Recall that a bounded operator on a Banach space is said to be {\it power compact} if some power is a compact operator.

\begin{proposition}
The zero operator is the only positive operator on a Hilbert lattice that is a self-commutator of a positive power compact operator.
\end{proposition}

\begin{proof}
Suppose that a positive operator $C$ on a Hilbert lattice $E$ is a self-commutator of a positive power compact operator $A$.
Since the operator $C=A^*A-AA^*$ is positive and $A$ is a positive power compact operator, by \cite[Lemma 2.3]{DK11} the operator $C$ is also power compact.

We claim that $C$ is a quasinilpotent operator. First, consider the complexifications $A_{\mathbb C}$ and $A^*_{\mathbb C}$ of $A$ and $A^*$, respectively. Since $A_{\mathbb C}$ and $A^*_{\mathbb C}$ are positive power compact operators that satisfy
$A^*_{\mathbb C}A_{\mathbb C}-A_{\mathbb C}A^*_{\mathbb C}\geq 0$, by \cite[Proposition 4.7]{KS17}, the operators $A^*_{\mathbb C}$ and $A_{\mathbb C}$ are simultaneously triangularizable. Suppose that $C^n$ is compact for some $n\in \mathbb N$. Then the operator
$$C_{\mathbb C}^n=(A^*_{\mathbb C}A_{\mathbb C}-A_{\mathbb C}A^*_{\mathbb C})^n$$ is compact on $E_{\mathbb C}$. Since $A_{\mathbb C}$ and $A^*_{\mathbb C}$ are simultaneusly triangularizable, the diagonal coefficients of the operators $C_{\mathbb C}^n$ are all zero. By Ringrose's theorem (see \cite[Theorem 1]{Ringrose} and \cite[Theorem 7.2.3]{RaRo}) the operator $C_{\mathbb C}^n$ is quasinilpotent. Hence, $C_{\mathbb C}$ and $C$ are quasinilpotent. Since $C$ is self-adjoint, it is the zero operator.
\end{proof}

The compactness assumption cannot be omitted in the last proposition, as the example of the unilateral shift operator on $\ell^2$ shows.

\section{Positive linear isometries between separable $L^p$-spaces}\label{Positive isometries}

In \Cref{Positive central self-commutators} we will consider the question which positive operators on Hilbert lattices are self-commutators of positive operators.
As a first result, we will prove that every positive compact central operator on a separable infinite-dimensional Hilbert lattice is a self-commutator of a positive operator. In our proof we will use the fact that for every positive diagonal operator $D$ on $\ell^2_n$ there exists a positive operator $X\colon \ell^2_n\to L^2[0,1]$ such that $D=X^*X$ (see \Cref{positive_square_root_inclusion_n}). A more general result holds also when one replaces $\ell^2_n$ by $\ell^2$.
As a special case, when $D=I$, we obtain that there always exist a positive isometry from $\ell^2\to L^2[0,1]$.

\begin{proposition}\label{positive_square_root_inclusion}
For every positive bounded diagonal operator $D$ on $\ell^2$ there exists a positive linear operator $X\colon \ell^2 \to L^2[0,1]$ such that
$X^*X=D$. Furthermore, there exists a positive self-adjoint operator $Y\colon L^2[0,1]\to L^2[0,1]$ such that $Y^2=XX^*$.
\end{proposition}

\begin{proof}
Let $(d_n)_{n\in \mathbb N}$ be the sequence of diagonal entries of the operator $D$ and let
$(e_n)_{n\in \mathbb N}$ be the standard basis of $\ell^2$.
We claim that $X\colon \ell^2 \to L^2[0,1]$ defined as
$$X\colon x\mapsto \sum_{n=1}^\infty \sqrt{2^nd_n}\langle x,e_n\rangle \chi_{n} , $$
where $\chi_n$ is the characteristic function of the interval $[\frac{1}{2^n},\frac{1}{2^{n-1}}]$ is a bounded operator which satisfies $X^*X=D$.
To see this, first define $M:=\sup_{n\in \mathbb N}d_n$ and note that
\begin{align*}
\|Xx\|^2&=\sum_{n=1}^\infty 2^nd_n|\langle x,e_n\rangle|^2 \int_{2^{-n}}^{2^{-n+1}}\hspace{-1em}dt=\sum_{n=1}^\infty 2^nd_n|\langle x,e_n\rangle|^2 2^{-n}\\
&\leq M\sum_{n=1}^\infty |\langle x,e_n\rangle|^2=M\|x\|^2
\end{align*}
shows that $X$ is bounded with $\|X\|\leq \sqrt M$.
From
$\langle X^*f,e_k\rangle=\langle f,Xe_k\rangle = \sqrt{2^k d_k} \langle f,\chi_k \rangle $  it follows
$$X^*f=\sum_{k=1}^\infty \sqrt{2^k d_k} \langle f,\chi_k \rangle e_k$$ and so
\begin{align*}
X^*Xx&=\sum_{n=1}^\infty \sqrt{2^n d_n}\langle x,e_n\rangle\, X^*\chi_n=\sum_{n=1}^\infty \sqrt{2^n d_n}\langle x,e_n\rangle\,\sum_{k=1}^\infty \sqrt{2^k d_k} \langle \chi_n,\chi_k \rangle e_k\\
&=\sum_{n=1}^\infty 2^nd_n \langle x,e_n\rangle 2^{-n} e_n=\sum_{n=1}^\infty d_n \langle x,e_n\rangle e_n\\
&=Dx.
\end{align*}
Similarly as above, a direct calculation shows that for every $f\in L^2[0,1]$ we have
$$XX^*f=\sum_{k=1}^\infty 2^kd_k\langle f,\chi_k\rangle \chi_k.$$
Hence, the positive operator $Y\colon L^2[0,1]\to L^2[0,1]$ defined as
$$Yf=\sum_{k=1}^\infty 2^k\sqrt{d_k}\langle f,\chi_k\rangle \chi_k$$
satisfies $Y^2=XX^*$ and $Y^*=Y$.
\end{proof}

If one replaces the Hilbert lattice $\ell^2$ with the finite-dimensional Hilbert lattice $\ell^2_n$, then the operator $X\colon x\mapsto \sum_{k=1}^n \langle x,e_k\rangle \sqrt{2^kd_k}\, \chi_{k}$ satisfies $X^*X=D$. Moreover, as above, a direct calculation shows that for every $f\in L^2[0,1]$ we have
$$XX^*f=\sum_{k=1}^n 2^kd_k\langle f,\chi_k\rangle \chi_k.$$
The positive operator $Y\colon L^2[0,1]\to L^2[0,1]$ defined as
$$Yf=\sum_{k=1}^n 2^k\sqrt{d_k}\langle f,\chi_k\rangle \chi_k$$
satisfies $Y^2=XX^*$ and $Y^*=Y$.
This finite-dimensional adaptation of the proof of \Cref{positive_square_root_inclusion} gives the following result.

\begin{corollary}\label{positive_square_root_inclusion_n}
For every positive $n\times n$ diagonal matrix $D$ there exists a positive linear operator $X\colon \ell^2_n\to L^2[0,1]$ such that $X^*X=D$. Furthermore, there exists a positive self-adjoint operator $Y\colon L^2[0,1]\to L^2[0,1]$ such that $Y^2=XX^*$.
\end{corollary}

\begin{corollary}\label{positive isometry existence}
There exist positive linear isometries $U\colon \ell^2\to L^2[0,1]$ and $U_n\colon \ell^2_n\to L^2[0,1]$ for each $n\in \mathbb N$.
\end{corollary}

\begin{proof}
By \Cref{positive_square_root_inclusion} there exists a positive linear operator $X\colon \ell^2\to L^2[0,1]$ such that $X^*X$ is the identity operator on $\ell^2$. Therefore, $U:=X\colon \ell^2\to L^2[0,1]$ is the desired positive linear isometry.
To find a positive isometry $U_n\colon \ell^2_n\to L^2[0,1]$ one can apply \Cref{positive_square_root_inclusion_n} instead of \Cref{positive_square_root_inclusion}.
\end{proof}

Since $\ell^2$ and $L^2[0,1]$ are isometrically isomorphic as Hilbert spaces,
a natural question which arises here is whether there exists a positive surjective linear isometry, i.e., a positive unitary operator $U\colon \ell^2\to L^2[0,1]$. Flores proved in \cite{Flores} that there even does not exist a regular isomorphism between $\ell^2$ and $L^2[0,1]$. The following proposition shows that there is no positive linear isometry from $L^2[0,1]$ to $\ell^2$. The obstruction comes from lattice structures of both Banach lattices.

\begin{proposition}\label{no_positive_isometries}
There is no positive linear isometry $V\colon L^{p}[0,1]\to \ell^{p}$ for any $1<p<\infty$.
\end{proposition}

\begin{proof}
Suppose there exists a positive linear isometry $V\colon L^p[0,1]\to \ell^p$. We claim that $V$ is a lattice homomorphism.
To see this, it suffices to prove that $f\wedge g=0$ in $L^p[0,1]$ implies $Vf\wedge Vg=0$ in $\ell^p$. Therefore, assume that
$f\wedge g=0$ in $L^p[0,1]$. Then
$$\|f+g\|_p^p=\int_0^1 (f+g)^p\,dt=\int_0^1 f^p\,dx+\int_0^1g^p\,dt=\|f\|_p^p +\|g\|_p^p$$ and so
$$\|V(f+g)\|_p^p=\|f+g\|_p^p=\|f\|_p^p+\|g\|_p^p=\|Vf\|_p^p+\|Vf\|_p^p.$$
Positivity of $Vf$ and $Vg$ implies
$$\sum_{n=1}^\infty \left((Vf)_n+(Vg)_n\right)^p=\sum_{n=1}^\infty (Vf)_n^p+\sum_{n=1}^\infty (Vg)_n^p.$$

We claim that, for each $n$,  one of $(Vf)_n$ and $(Vg)_n$ is nonzero only if the other one is zero. To prove this, it suffices to prove that
for positive real numbers $a,b$ and $1< p<\infty$ we always have $(a+b)^p>a^p+b^p$. Let us introduce the function
$f\colon [0,\infty)\to \mathbb R$ defined as $f(x):=(x+b)^p-x^p-b^p$. Then $f'(x)=p\left((x+b)^{p-1}-x^{p-1}\right)>0$. Hence, $f$ is strictly increasing, and so $(a+b)^p-a^p-b^p=f(a)>f(0)=0$. This shows that, for each $n$, $(Vf)_n$ and $(Vg)_n$ cannot be positive simultaneously. This proves that $Vf\wedge Vg=0$, and so, $V$ is a lattice homomorphism.

Since $V$ is an isometric lattice homomorphism, the Banach lattice $\ell^p$ contains a closed sublattice $E$ which is lattice isometric to the Banach lattice $L^p[0,1]$. As such, $E$ needs to be without atoms. On the other hand, \cite[Theorem 3.2]{TroitskyBilokopytov} implies that $E$ is the closed span of an infinite disjoint positive sequence, and so $E$ is atomic. This contradiction shows that $V$ cannot exist.
\end{proof}

It should be noted that the proof of \Cref{no_positive_isometries} shows that there is no positive operator $T\colon L^p[0,1]\to \ell^p$ such that the restriction of $T$ to the positive cone $L^p[0,1]_+$ is an isometry. For $p=1$ such operators exist.

\begin{example}
Let us define the linear operator $V\colon L^1[0,1]\to \ell^1$ by
$$(Vf)_n=\int_{2^{-n}}^{2^{-n+1}}f(t)\,dt.$$
Since
$$\sum_{n=1}^\infty |(Vf)|_n\leq \sum_{n=1}^\infty\int_{2^{-n}}^{2^{-n+1}}|f(t)|\,dt=\int_{0}^1 |f(t)|\,dt=\|f\|_1,$$
it follows that $\|Vf\|_1\leq \|f\|_1$, and so, $V$ is a contraction. Moreover, if $f$ is a non-negative function, then
the calculation above yields $\|Vf\|_1=\|f\|_1$. Positivity of $V$ should be clear.
\end{example}

Although \Cref{no_positive_isometries} shows that there does not exist a positive linear isometry $V\colon L^p[0,1]\to \ell^p$ for $1<p<\infty$, spaces $\ell^2$ and $L^2[0,1]$ are isometrically isomorphic as Hilbert spaces. For $p=1$ the situation is more severe as there is no linear isometry from $L^1[0,1]$
to $\ell^1$. This follows from the fact that by \cite[Example 3.1.2]{Albiac} the Banach space $\ell^1$  has an unconditional basis, so that by \cite[Theorem 6.3.3]{Albiac} the Banach space $L^1[0,1]$ cannot be embedded into a Banach space with an unconditional basis.

\section{Positive central self-commutators}\label{Positive central self-commutators}

In this section we consider the question which positive central operators on separable infinite-dimensional Hilbert lattices are positive self-commutators of positive operators. We start by proving an order analog of \cite[Theorem 5]{Radjavi:66}. 

\begin{theorem}\label{compact operators are self-commutators}
Every positive compact central operator on a separable infinite-dimensional Hilbert lattice is a self-commutator of a positive operator.
\end{theorem}

\begin{proof}
Let $C$ be a positive compact central operator on a separable infinite-dimensional Hilbert lattice $\mathcal H$.
As mentioned in \Cref{Introduction},
%By Bohnenblust's result (\cite[Theorem 7.1]{Bohnenblust}),
$\mathcal H$ is lattice isometric to one of the following Hilbert lattices:
$\ell^2$, $L^2[0,1]$, $\ell^2\oplus L^2[0,1]$ and $\ell^2_n\oplus L^2[0,1]$. We will consider each of these cases separately.

\emph{Case} $\mathcal H\cong \ell^2$: Since $C$ is a central operator on $\ell^2$, it is a diagonal operator with positive diagonal coefficients $(d_n)_{n\in\mathbb N}$. From compactness of $C$ it follows that $d_n\to 0$ as $n$ goes to infinity, so that there exists the largest diagonal coefficient. Without loss of generality we may assume that $d_1\geq d_k$ for every $k\in \mathbb N$.  Since for every $\epsilon>0$ there are only finitely many $n\in \mathbb N$ such that $d_n\geq \epsilon$  there exist countably infinitely many sequences $s^{(n)}$  such that
\begin{itemize}
\item[(i)] the terms of all sequences $s^{(n)}$ equal the terms of the sequence $(d_n)_{n\in\mathbb N}$ with respect to all of their multiplicities;
\item[(ii)] the first term $s_1^{(n)}$ of the sequence $s^{(n)}$ is its largest term for each $n\in \mathbb N$;
\item[(iii)] $s_1^{(n)}\leq \frac{d_1}{2^{n-1}}$ for every $n\in \mathbb N$.
\end{itemize}

Let $e_k^{(n)}$ be the standard basis eigenvector that corresponds to the eigenvalue $s_k^{(n)}$ for $n,k\in \mathbb N$ and let $\mathcal H_n$ be the closed linear span of $\{e_k^{(n)}:\; k\in \mathbb N\}$. Then every $\mathcal H_n$ is lattice isometric to $\ell^2$. Let us define the isometric lattice isomorphism   $U\colon \ell^2 \to \left(\bigoplus_{n=1}^\infty \mathcal H_n\right)_2$ as follows. For each $m\in \mathbb N$ there exist unique positive integers $k,n\in \mathbb N$ such that $e_m=e_k^{(n)}$.  We define $Ue_m=e_k^{(n)}$ and extend $U$ by linearity on the linear span of the set of all standard basis vectors. Clearly, the extended mapping is a positive bounded linear operator which can be uniquely extended by continuity to a positive operator on  $\ell^2$. Since by construction $U$ maps the set $\{e_n:\; n\in \mathbb N\}$ injectively onto itself, we have that $U$ is unitary. Moreover, an easy calculation shows that
$$UCU^{-1}=
\left( \begin{matrix}
D_1 & 0 & 0  & \hdots\\
0 & D_2 & 0  & \hdots\\
0 & 0 & D_3  & \hdots\\
\vdots & \vdots & \vdots & \ddots
\end{matrix}
\right)$$
where $D_n$ is the diagonal operator $\diag\big(\{s_k^{(n)}:\;k\in \mathbb N\}\big)$ on $\mathcal H_n$.

Let us solve an infinite system of operator equations given by $X_1^2=D_1$ and $X_{n+1}^2-X_n^2=D_{n+1}$ for $n\in \mathbb N$.
If $n=1$, the positive diagonal operator $X_1:=\sqrt{D_1}$ satisfies $X_1^2=D_1$. If $n=2$, then the positive diagonal operator $X_2:=\sqrt{D_1+D_2}$ solves the operator equation $X_2^2=X_1^2+D_2$. Inductively, one can show that for each $n\in \mathbb N$ the positive diagonal operator $X_n:=\sqrt{D_1+\cdots+D_n}$ formally solves the operator equation $X_n^2-X_{n-1}^2=D_n$.
Since the norm of the operator $X_n$ satisfies
$$\|X_n\|=\sqrt{s_1^{(1)}+s_1^{(2)}+\cdots+s_1^{(n)}}\leq \sqrt{d_1+\frac{d_1}{2}+\cdots +\frac{d_1}{2^{n-1}}}\leq \sqrt{2d_1},$$
the operator $X_n$ is bounded.

Let us define the operator $X$ on $\ell^2\cong \big(\bigoplus_{n=1}^\infty \ell^2\big)_2$ as the infinite operator matrix
$$X=\left( \begin{matrix}
0 & 0 & 0 &  \hdots\\
X_1 & 0 & 0 &  \hdots\\
%0 &  0 & X_3 & \hdots\\
0 &  X_2 & 0 & \hdots\\
\vdots & \vdots & \vdots & \ddots
\end{matrix}
\right).$$
Since $\|X_n\|\leq \sqrt{2d_1}$ for each $n\in \mathbb N$, the operator $X$ is bounded with $\|X\|\leq \sqrt{2d_1}$.
Self-adjointness of $X_n$ now yields
$$[X^*,X]=\left( \begin{matrix}
X_1^2  & 0 & 0  & \hdots\\
0 & X_2^2-X_1^2 & 0  & \hdots\\
0 & 0 & X_3^2-X_2^2  & \hdots\\
\vdots & \vdots & \vdots & \ddots
\end{matrix}
\right)=\left( \begin{matrix}
D_1 & 0 & 0  & \hdots\\
0 & D_2 & 0  & \hdots\\
0 & 0 & D_3  & \hdots\\
\vdots & \vdots & \vdots & \ddots
\end{matrix}
\right),$$
and so $UCU^{-1}=[X^*,X]$. Since $U$ satisfies $U^*=U^{-1}$, we finally have $C=U^{-1}[X^*,X]U=[U^{-1}X^*U,U^{-1}XU]=[(U^{-1}XU)^*,U^{-1}XU]$.

\emph{Case} $\mathcal H\cong L^2[0,1]$: Since the Banach lattice $L^2[0,1]$ is atomless, by \cite[Corollary 1.7]{Schep80} it follows that $C=0$ on $L^2[0,1]$, so that it is a self-commutator of itself.

\emph{Case} $\mathcal H\cong \ell^2\oplus L^2[0,1]$: Since $C$ is a positive central operator on $\ell^2\oplus L^2[0,1]$, it can be written as $C=C_1\oplus C_2$ where $C_1$ and $C_2$ are positive central operators on $\ell^2$ and $L^2[0,1]$, respectively. Compactness of $C$ yields compactness of $C_1$ and $C_2$. The previous two cases imply that $C_2=0$ and $C_1=[X^*,X]$ for some positive operator $X$ on $\ell^2$. Therefore, the operator $C=[X^*,X]\oplus 0=[X^*\oplus 0,X\oplus 0]$ is a self-commutator.

\emph{Case} $\mathcal H\cong \ell^2_n\oplus L^2[0,1]$: Since $C$ is a positive central operator on $\ell^2_n\oplus L^2[0,1]$, it can be written as $C=C_1\oplus C_2$ where $C_1$ and $C_2$ are positive central operators on $\ell^2_n$ and $L^2[0,1]$, respectively. Similarly, as in the previous case, it follows that $C_2=0$, and so $C=C_1\oplus 0$ for some positive diagonal $n\times n$ matrix $C_1$.

By \Cref{positive_square_root_inclusion_n} there exist a positive operator $X\colon \ell^2_n\to L^2[0,1]$ such that $C_1=X^*X$ and a positive self-adjoint operator $Y\colon L^2[0,1]\to L^2[0,1]$ such $XX^*=Y^2$. With respect to the decomposition $\mathcal H\cong \ell^2_n\oplus \bigoplus_{n=1}^\infty L^2[0,1]$ we define the positive operator $Z$ as an infinite block operator matrix
$$Z=\left( \begin{matrix}
0 & 0 & 0 & 0 & \hdots\\
X & 0 & 0 & 0 & \hdots\\
0 & Y & 0 & 0 & \hdots\\
0 & 0 & Y & 0 & \hdots\\
\vdots & \vdots & \vdots & \vdots & \ddots
\end{matrix}
\right).$$
A direct calculation  shows
$$Z^*Z=
\left( \begin{matrix}
X^*X & 0 & 0 & 0 & \hdots\\
0 & Y^2  & 0 & 0 & \hdots\\
0 & 0 & Y^2  & 0 & \hdots\\
0 & 0 & 0 &  Y^2 &  \hdots\\
\vdots & \vdots & \vdots & \vdots & \ddots
\end{matrix}
\right) \qquad \textrm{and} \qquad ZZ^*=
\left( \begin{matrix}
0 & 0 & 0 & 0 & \hdots\\
0 & XX^* & 0 & 0 & \hdots\\
0 & 0 & Y^2 & 0 & \hdots\\
0 & 0 &0 & Y^2 & \hdots\\
\vdots & \vdots & \vdots & \vdots & \ddots
\end{matrix}
\right).$$
Since $X^*X=C_1$ and $Y^2=XX^*$, we have
\[
[Z^*,Z]=Z^*Z-ZZ^*=\left( \begin{matrix}
C_1 & 0 & 0 & 0 & \hdots\\
0 & 0 & 0 & 0 & \hdots\\
0 & 0 & 0  & 0 & \hdots\\
0 & 0 & 0 &  0 &  \hdots\\
\vdots & \vdots & \vdots & \vdots & \ddots
\end{matrix}
\right)=C
\]%
which finishes the proof.
\end{proof}

\begin{corollary}
Every positive compact central operator of infinite rank on a Hilbert lattice is a self-commutator of a positive operator.
\end{corollary}

\begin{proof}
Let $C$ be a central operator on a Hilbert lattice $\mathcal H$. The range $\ran C$ of a central operator $C$ is an ideal in $\mathcal H$, so that it is equal to the range ideal $\cR(C)$. Since the closure $\overline{\cR(C)}$ is a closed ideal in an order continuous Banach lattice $\mathcal H$, we have that $\overline{\cR(C)}$ is a band in $\mathcal H$. The band generated by $\cR(C)$ is $\cR(C)^{dd}$, so that $\cR(C)^{dd}\subseteq \overline{\cR(C)}$. On the other hand, every band in a normed lattice is a closed ideal from where it follows $\overline{\cR(C)} \subseteq \cR(C)^{dd}$. Since every band in $\mathcal H$ is a projection band, we have
$$\mathcal H=\cR(C)^{dd}\oplus \cR(C)^{d}=\overline{\cR(C)}\oplus \cR(C)^{d}.$$
With respect to this decomposition, $C$ can be written as $C=C_1\oplus 0$ where $C_1=C|_{\overline{\cR(C)}}$. Since $C_1$ is a compact operator, its range is separable, so that $\overline{\cR(C)}$ is infinite-dimensional and separable as well. Hence, by \Cref{compact operators are self-commutators} the operator $C_1$ is a self-commutator of a positive operator. To finish the proof we observe that the zero operator is self-commutator of itself.
\end{proof}

Our next goal is to prove that positive central operators with infinite-dimensional kernels on a Hilbert lattice which is either $\ell^2$ or $L^2[0,1]$ are always self-commutators of positive operators. Before we prove this result, we need the following lemma.

\begin{lemma}\label{o pozitivnem sebikomutatorju}
Let $T$ be a positive operator on a Hilbert lattice $\mathcal H$ which is either $\ell^2$ or $L^2[0,1]$. If $T$ is positive semi-definite operator whose square root is a positive operator, then the operators
$$\left( \begin{matrix}
T & 0 \\
0 & 0
\end{matrix}
\right) \qquad \textrm{and}\qquad
\left( \begin{matrix}
0 & 0\\
0 & T
\end{matrix}
\right)$$  on the Hilbert lattice $\mathcal H\oplus \mathcal H$
are self-commutators of positive operators.
\end{lemma}

\begin{proof}
We only consider the case $\mathcal H=\ell^2$ as the case of $\mathcal H=L^2[0,1]$ can be treated similarly.
Let $\sqrt T$ be the square-root of $T$.
We define the operator $A$ on the Hilbert lattice $\ell^2\oplus \ell^2\cong \ell^2\oplus \bigoplus_{n=1}^\infty \ell^2$ as
$$A=\left( \begin{matrix}
0 & 0 & 0 & 0 & \hdots\\
\sqrt{T} & 0 &0 & 0 & \hdots\\
0 & \sqrt{T} & 0 & 0 & \hdots\\
0 & 0 & \sqrt{T} & 0 & \hdots\\
\vdots & \vdots & \vdots & \vdots & \ddots
\end{matrix}
\right).$$
Since $\sqrt T$ is self-adjoint, a direct calculation shows that $[A^*,A]=T\oplus 0$.
Positivity of $\sqrt T$ yields that $A$ is a positive operator on $\ell^2\oplus \ell^2$.

To prove that the operator $0\oplus T$ is a self-commutator of a positive operator, we first introduce the switch operator $U\colon \mathcal H\oplus \mathcal H\to \mathcal H\oplus \mathcal H$ given by the $2\times 2$ block-operator matrix
$$\left( \begin{matrix}
0 & I\\
I & 0
\end{matrix}
\right).$$
Clearly, $U$ is an isometric self-adjoint isomorphism of the Hilbert lattice $\mathcal H\oplus \mathcal H$ which satisfies $U=U^{-1}=U^*$.
Since
$$U^*\left( \begin{matrix}
0 & 0 \\
0 & T
\end{matrix}
\right)U=\left( \begin{matrix}
T & 0 \\
0 & 0
\end{matrix}
\right),$$ by the case above, there exists a positive operator $A$ on $\mathcal H\oplus \mathcal H$ such that
$T\oplus 0=[A^*,A]$, and so
$0\oplus T=U(T\oplus 0)U^*=U[A^*,A]U^*=[UA^*U^*,UAU^*]$.
\end{proof}

We now prove an order analog of \cite[Theorem 3]{Radjavi:66}.

\begin{theorem}\label{SelfCommutatorsL^2}
Let $C$ be a positive central operator on a separable infinite-dimensional Hilbert lattice $\ell^2$ or $L^2[0,1]$ with infinite dimensional kernel. Then $C$ is a self-commutator of a positive operator.
\end{theorem}

\begin{proof}
Obviously, we may assume that $C\neq 0$ as the zero operator is a self-commutator of itself.
Since $C$ is a central operator, it is order continuous and its kernel $\ker C$ equals its absolute kernel which is a projection band.
Its orthogonal complement $(\ker C)^\perp =(\ker C)^d$ is the band generated by the range of $C$. Since $\ell^2=(\ker C)^\perp\oplus \ker C$, the restriction of $C$ to $(\ker C)^\perp$ is a strictly positive central operator.

Suppose that the underlying Hilbert lattice is $\ell^2$. If $\dim (\ker C)^\perp$ is finite, then $C$ is a compact operator, so that it is a self-commutator of a positive operator by \Cref{compact operators are self-commutators}. If $\dim (\ker C)^\perp$ is infinite, then we can decompose $\ell^2=(\ker C)^\perp \oplus  \ker C\cong  \ell^2\oplus \ell^2$ and the block-operator matrix corresponding to $C$ is of the form
$$\left( \begin{matrix}
C|_{(\ker C)^\perp} & 0 \\
0 & 0
\end{matrix}
\right).$$
% Since $C$ is a central operator, $\ker C$ is a band in $\ell^2$, so that $(\ker C)^\perp=(\ker C)^d$.
Here the decomposition $\ell^2=(\ker C)^\perp \oplus  \ker C\cong  \ell^2\oplus \ell^2$ is simultaneously a Hilbert space and a Banach lattice decomposition. This implies that $C|_{(\ker C)^\perp}$ is a diagonal operator with positive diagonal entries on the Hilbert lattice $(\ker C)^\perp\cong \ell^2$. As such, it is positive semi-definite operator which admits a positive square root.
\Cref{o pozitivnem sebikomutatorju} yields that $C$ is a self-commutator of a positive operator.

When the underlying Hilbert lattice is $L^2[0,1]$ we argue similarly.
Since $C\neq 0$, we can decompose $L^2[0,1]=(\ker C)^\perp \oplus \ker C\cong L^2[0,1]\oplus L^2[0,1]$. Then, with respect to this decomposition, the operator $C$ can be represented as a $2\times 2$ block-operator matrix
$$\left( \begin{matrix}
C|_{(\ker C)^\perp} & 0 \\
0 & 0
\end{matrix}
\right).$$
Since the operator $C|_{(\ker C)^\perp}$ is a positive central operator on $L^2[0,1]$, it is a multiplication operator $M_\varphi$ for some non-negative function $\varphi \in L^\infty[0,1]$. Since $M_\varphi$ is positive semi-definite with the positive square root $M_{\sqrt{\varphi}}$, another application of \Cref{o pozitivnem sebikomutatorju} yields that $C$ is a self-commutator of a positive operator.
\end{proof}

\begin{corollary}\label{band projections self-commutators}
A band projection on $\ell^2$ or $L^2[0,1]$ is a self-commutator of a positive operator if and only if its kernel is infinite-dimensional.
\end{corollary}

\begin{proof}
If $P$ is a band projection with a finite-dimensional kernel, then $I-P$ is a finite-rank projection. Hence, by an application of Wintner-Wielandt's result for the Calkin algebra (see \cite{Ha82}), the operator $P=I-(I-P)$ is not a commutator. The other implication follows from \Cref{SelfCommutatorsL^2}.
\end{proof}

Since $L^p(\mu)$ is not a Hilbert space when $p\neq 2$ we cannot talk about self-commutators. Nevertheless, the proof of  \Cref{SelfCommutatorsL^2} reveals the following result.

\begin{proposition}
Let $C$ be a positive central operator on a separable infinite-dimensional Banach lattice $\ell^p$ or $L^p[0,1]$ ($1 \le p < \infty$). If $\ker C$ is infinite dimensional, then $C$ is a commutator of positive operators.
\end{proposition}

We conclude this section with a result stating that invertible positive central operators on Banach lattices are not commutators of positive operators.
In the proof of this result we will need the following slight improvement of \cite[Theorem 2.1]{DK25} which for the sake of simplicity we state only for the case of operators on Banach lattices. 

\begin{theorem}\label{dominira identiteto}
Let $A$ and $B$ be operators on a Banach lattice $L$ with one of them positive or negative. Then there does not exist a positive scalar $c>0$ such that
$AB-BA\geq cI$.
\end{theorem}

%This result was proved by the authors in \cite[Theorem 2.1]{DK25} for $c=1$. The general case follows by homogeneity.

\begin{theorem}
An invertible positive central operator on a Banach lattice is not a commutator of two operators with one of them positive or negative.
\end{theorem}

\begin{proof}
Since the center $\mathscr Z(L)$ of a Banach lattice $L$ is isomorphic as an $f$-algebra to an $f$-algebra $C(K)$ for some compact Hausdorff space $K$, a positive central operator $C$ is invertible if and only if there exists a positive scalar $\delta>0$ such that $C\geq \delta I$. Hence, \Cref{dominira identiteto} yields that there do not exist bounded operators $A$ and $B$ on $L$ with one of them positive or negative such that $C=AB-BA$.
\end{proof}

\section{Positive central operators are sums of two positive self-commutators}\label{Section 5}

It was proved in \Cref{SelfCommutatorsL^2} that every positive central operator with infinite-dimensional kernel on a Hilbert lattice $\ell^2$ or $L^2[0,1]$ is a self-commutator of a positive operator. It turns out that the situation becomes very difficult when we consider the separable Hilbert lattice $\ell^2\oplus L^2[0,1]$. To explain where the difficulty lies, let us consider the following example.

\begin{example}\label{example C}
The operator $C$ on a Hilbert lattice $\mathcal H=\ell^2\oplus L^2[0,1]$ represented with the $2\times 2$ operator matrix
$$\left(\begin{matrix}
I & 0\\
0 & 0
\end{matrix}\right)$$
is a self-commutator of a positive operator on $\mathcal H$.
To see this, observe first that by \Cref{positive_square_root_inclusion} there exists a positive isometry $X\colon \ell^2\to L^2[0,1]$ and a positive self-adjoint operator $Y\colon L^2[0,1]\to L^2[0,1]$ such that $Y^2=XX^*$.
With respect to the Hilbert lattice decomposition $\mathcal H\cong \ell^2\oplus \bigoplus_{n=1}^\infty L^2[0,1]$, as in \Cref{compact operators are self-commutators} we define the positive operator $Z$ as the infinite block operator matrix
$$Z=\left( \begin{matrix}
0 & 0 & 0 & 0 & \hdots\\
X & 0 & 0 & 0 & \hdots\\
0 & Y & 0 & 0 & \hdots\\
0 & 0 & Y & 0 & \hdots\\
\vdots & \vdots & \vdots & \vdots & \ddots
\end{matrix}
\right).$$
An easy calculation shows that $[Z^*,Z]=C$ proving that $C$ is a self-commutator of a positive operator.
\end{example}

At first glance it may seem that a similar idea works to prove that the operator $C$ on $\mathcal H=\ell^2\oplus L^2[0,1]$ represented with a $2\times 2$ block operator matrix
$$\left(\begin{matrix}
0 & 0\\
0 & I
\end{matrix}\right)$$
is a self-commutator of a positive operator. However, there is a big issue. The proof that the operator from \Cref{example C} is a self-commutator of a positive operator heavily relies on the fact that there exists a positive isometry $X\colon \ell^2\to L^2[0,1]$. However, by \Cref{no_positive_isometries} there are no positive isometries from $L^2[0,1]$ to $\ell^2$.

In this section we prove that every positive central operator on a separable infinite-dimensional Hilbert lattice is a sum of two positive self-commutators of positive operators.

\begin{theorem}
Every positive central operator $C$ on an infinite-dimensional separable Hilbert lattice is a sum of two positive self-commutators of positive operators.
\end{theorem}

\begin{proof}
Suppose first that $\mathcal H$ is either $\ell^2$ or $L^2[0,1]$.
Then, $\mathcal H\cong \mathcal H\oplus \mathcal H$ as a Hilbert lattice. Since $C$ is a central operator, with respect to the decomposition $\mathcal H\oplus \mathcal H$ the operator $C$ can be represented as a $2\times 2$ operator matrix
$$\left( \begin{matrix}
C_1 & 0 \\
0 & C_2
\end{matrix}
\right)=\left( \begin{matrix}
C_1 & 0 \\
0 & 0
\end{matrix}
\right)+\left( \begin{matrix}
0 & 0 \\
0 & C_2
\end{matrix}
\right)$$
where $C_1$ and $C_2$ are positive central operators on $\mathcal H$. If $\mathcal H=\ell^2$, then $C_1$ and $C_2$ are diagonal operators with positive diagonal entries. If $\mathcal H=L^2[0,1]$, then $C_1$ and $C_2$ are multiplication operators with almost everywhere non-negative multiplying functions. Since positive diagonal operators on $\ell^2$ and positive multiplication operators on $L^2[0,1]$ are positive and positive semi-definite whose roots are positive operators, by \Cref{o pozitivnem sebikomutatorju} there exist positive operators $A$ and $B$ on $\mathcal H\oplus \mathcal H$ such that
$$\left( \begin{matrix}
C_1 & 0 \\
0 & C_2
\end{matrix}
\right)=[A^*,A]+[B^*,B].$$

Now we consider the case when $\mathcal H=\ell^2\oplus L^2[0,1]$. Since $C$ is a positive central operator on $\mathcal H$, with respect to the decomposition $\ell^2\oplus L^2[0,1]$ we can write
$$C=\left( \begin{matrix}
C_1 & 0 \\
0 & C_2
\end{matrix}
\right)$$
where $C_1$ and $C_2$ are positive central operators on $\ell^2$ and $L^2[0,1]$, respectively. By the cases above, there exist positive operators $A_1,B_1$ on $\ell^2$ and $A_2,B_2$ on $L^2[0,1]$ such that $C_1=[A_1^*,A_1]+[B_1^*,B_1]$ and $C_2=[A_2^*,A_2]+[B_2^*,B_2]$.
Hence,
\begin{align*}
\left( \begin{matrix}
C_1 & 0 \\
0 & C_2
\end{matrix}
\right)&=
\left( \begin{matrix}
[A_1^*,A_1]+[B_1^*,B_1] & 0\\
0 & [A_2^*,A_2]+[B_2^*,B_2]
\end{matrix}
\right)\\
&=
\left( \begin{matrix}
[A_1^*,A_1] & 0 \\
0 & [A_2^*,A_2]
\end{matrix}
\right)+\left( \begin{matrix}
[B_1^*,B_1] & 0 \\
0 & [B_2^*,B_2]
\end{matrix}
\right).
\end{align*}
The last expression can be written as the sum of two self-commutators of positive operators $A_1\oplus A_2$ and $B_1\oplus B_2$.

Finally, we consider the last case when the underlying Hilbert lattice $\mathcal H$ is of the form $\ell_n^2\oplus L^2[0,1]$ for some $n\in \mathbb N$. As above, due to the Hilbert lattice isomorphism $L^2[0,1]\cong L^2[0,1]\oplus L^2[0,1]$ we can furthermore decompose
$\ell_n^2\oplus L^2[0,1]\cong \ell_n^2\oplus L^2[0,1]\oplus L^2[0,1]$ and with respect to this decomposition
we can write
$$C=\left( \begin{matrix}
D & 0 & 0\\
0 & M_\varphi & 0 \\
0 & 0 & M_\psi
\end{matrix}
\right)=\left( \begin{matrix}
D & 0 & 0\\
0 & M_\varphi & 0 \\
0 & 0 & 0
\end{matrix}
\right)+\left( \begin{matrix}
0 & 0 & 0\\
0 & 0 & 0 \\
0 & 0 & M_\psi
\end{matrix}
\right)$$
for some non-negative functions $\varphi$ and $\psi\in L^\infty[0,1]$ and some non-negative diagonal $n \times n$ matrix $D$.
Since an application of \Cref{o pozitivnem sebikomutatorju} yields that the last operator above is a positive self-commutator of a positive operator,  due to the Hilbert lattice isomorphism $\ell_n^2\oplus L^2[0,1]\oplus L^2[0,1] \cong \ell_n^2\oplus L^2[0,1]\oplus L^2[0,1]\oplus L^2[0,1]$ it suffices to prove that  the operator given by the block-operator matrix
$$ \left( \begin{matrix}
D & 0 & 0\\
0 & M_\varphi & 0 \\
0 & 0 & 0
\end{matrix}
\right) \cong
\left( \begin{matrix}
D & 0 & 0 & 0\\
0 & M_\varphi  & 0 & 0\\
0 & 0 & 0 & 0\\
0 & 0 & 0 & 0
\end{matrix}\right) $$
is a positive self-commutator.
Furthermore, by applying the switch operator we need to prove that the operator
$$\left( \begin{matrix}
D & 0 & 0 & 0\\
0 & 0 & 0 & 0\\
0 & 0 & M_\varphi & 0\\
0 & 0 & 0 & 0
\end{matrix}
\right)$$ is a positive self-commutator of a positive operator.
By \Cref{compact operators are self-commutators} and \Cref{SelfCommutatorsL^2}, respectively, there exist positive operators $A$ and $B$ such that
$$\left( \begin{matrix}
 D & 0\\
0 & 0
\end{matrix}
\right)=[A^*,A]\qquad \textrm{and}\qquad \left( \begin{matrix}
 M_\varphi & 0\\
0 & 0
\end{matrix}\right)=[B^*,B].$$
To conclude the proof, observe that the operator
$$\left( \begin{matrix}
D & 0 & 0 & 0\\
0 & 0 & 0 & 0\\
0 & 0 & M_\varphi & 0\\
0 & 0 & 0 & 0
\end{matrix}
\right)$$
is a self-commutator of the positive operator $A\oplus B$.
\end{proof}

We conclude this paper with the following question.

\begin{question}
Is the operator $C$ on the Hilbert lattice $\ell^2\oplus L^2[0,1]$ represented with a $2\times 2$ block operator matrix
$$\left(\begin{matrix}
0 & 0\\
0 & I
\end{matrix}\right)$$
a self-commutator of a positive operator?
\end{question}

\subsection*{Acknowledgements}
The first author was supported by the Slovenian Research and Innovation Agency program P1-0222.
The second author was supported by the Slovenian Research and Innovation Agency program P1-0222 and grant N1-0217.

\end{document}